% Logic Eprints
%Submitted 2213 Wed Sep 15, 1993 by: zapletal@math.psu.edu (jindrich zapletal )
%logic/zapletal/incns
%

\magnification=\magstep1
\input amstex
\documentstyle{amsppt}
\topmatter
\title
A New Proof of Kunen's Inconsistency
\endtitle
\author
Jind\v rich Zapletal
\endauthor
\affil
The Pennsylvania State University
\endaffil
\address
Department of Mathematics,
The Pennsylvania State University,
University Park, PA 16802
\endaddress
\email
zapletal\@math.psu.edu
\endemail
\abstract
Using elementary pcf, we show that there is no $j:V\to M,$ with 
$M$ transitive, $j\lambda=\lambda>crit (j),$
 $j^{\prime \prime }\lambda \in M.$
\endabstract
\endtopmatter

\document
For contradiction, assume there is such $j:V\to M.$ 
Let $\lambda$ be 
the least fixed point of $j$ above $crit(j).$ 
$j^{\prime \prime }\lambda \in M.$
It is easily established \cite {2} that $\lambda$ is a strong limit cardinal
of cofinality $\omega$ and $j(\lambda ^+)=\lambda ^+.$ We use the following 
fact due to Shelah.
\proclaim
{Lemma} If $\lambda >cof (\lambda )=\omega$ is a strong limit cardinal, 
there is a sequence of regular cardinals $\langle \lambda _i:i<\omega 
\rangle$ less than $\lambda ,$ with $sup ( \langle \lambda _i:i<\omega 
\rangle )=\lambda$ and $tcf \prod _{i<\omega}\lambda _i/fin =\lambda ^+.$
\endproclaim
For an elementary proof of this lemma see
 T. Jech's paper in \cite {1}. Fix $\langle 
\lambda _i:i<\omega \rangle$ converging to $\lambda$ as in the lemma.
W.l.o.g. $crit (j)<\lambda _i<\lambda$ for all $i<\omega .$ Fix a sequence
$G=\langle g_\alpha :\alpha <\lambda ^+\rangle \subset \prod _{i<\omega}
\lambda _i$ increasing and cofinal in $\prod _{i<\omega}\lambda _i/fin.$
By elementarity, $M\models$``$jG$ is increasing and cofinal in $j(\prod
_{i<\omega }\lambda _i/fin)".$ Define $g\in j(\prod _{i<\omega }\lambda _i)$
by $g(j\lambda _i)=sup(j^{\prime \prime}\lambda \cap j\lambda _i)
<\lambda _i.$ 
The last inequality follows from regularity of $j\lambda _i$
in $M.$ If $f\in \prod _{i<\omega} \lambda _i$ then $g>j(f)$ everywhere as
 $j(f)=j^{\prime \prime }f.$ Now $j^{\prime \prime }\lambda ^+$ is cofinal in 
$\lambda ^+$ since $\lambda ^+$ is fixed by $j$ and 
therefore 
$j^{\prime \prime }G$ is cofinal in $jG$ and thus in $j(\prod _{i<\omega }
\lambda _i/fin).$ However, we have just seen that $g\in j(\prod _{i<\omega }
\lambda _i)$ dominates every member of $j^{\prime \prime }G,$ contradiction.

\enddocument